\documentclass[11pt]{article}

\usepackage[all]{xy}

\usepackage{amsmath}

\usepackage{amssymb}

\usepackage{amscd}

\usepackage{bbm}

\newtheorem{theorem}{Theorem}[section]
\newtheorem{lemma}[theorem]{Lemma}

\newtheorem{exa}[theorem]{Example}

\newtheorem{exas}[theorem]{Examples}

\newtheorem{prope}[theorem]{Property}

\newtheorem{defini}[theorem]{Definition}

\newtheorem{rema}[theorem]{Remark}

\newenvironment{equationth}{\stepcounter{theorem}\begin{equation}}{\end{equation}}

\newcommand{\qed}{\ensuremath{\Box}}

\def\ind{{\rm{Ind}}}

\def\Z{\mathbbm{Z}}

\def\R{\mathbbm{R}}

\def\sph{\mathbb{S}}

\def\D{\mathbbm{D}}

\def\C{\mathbbm{C}}

\def\T{\mathbbm{T}}

\def\B{\mathbbm{B}}

\def\F{\mathcal F}

\def\a{\alpha}

\def\b{\beta }

\def\ro{\varrho}

\def\e{\varepsilon}

\def\D{\mathcal D}

\def\v {\vskip.1cm \noindent}

\def\vv {\vskip.2cm \noindent}

\def\ind{{Ind_{\rm PH}}}

\makeindex

\begin{document}

\title { Remarks on contact structures and vector fields
 on isolated complete intersection  singularities}

%%\endtitle

\author{  {Jos\'e Seade}\thanks{Supported by CONACYT and DGAPA-UNAM, Mexico)} \\  \\
 }

%\date{}

\date { }

\setcounter{section}{-1}

%%\NoRunningHeads

%\pagestyle{headings}

%\pagenumbering{roman}

\maketitle

\begin{abstract}
Let $(X,0)$ be an isolated complete intersection complex singularity ($X$ can also
be smooth at $0$). Let $K$ be 
its link, $\cal X$ its canonical contact structure and $\D_X$ the complex vector
bundle associated to $\cal X$. We prove that the bundle $\D_X$ is trivial  
 if and only if the Milnor number of $X$ satisfies $\mu(X,0) \equiv (-1)^{n-1}$
modulo $(n-1)!$. This follows from  a general theorem stating that the 
complex orthogonal complement of a vector field in $X$ with an isolated singularity
at $0$ is trivial iff the GSV-index of $v$ is a multiple of $(n-1)!$. 
We have also  an application to foliation theory: a holomorphic foliation
$\cal F$ in a ball $\B_r$ around the origin in $\C^3$, with an isolated 
singularity at $0$, admits a $C^\infty$ normal section (away from $0$) iff 
its multiplicity (or local index) is 
even, and this happens iff its normal bundle in $\B_r \setminus \{0\}$ is 
topologically trivial. 
\end{abstract}

\section{Introduction}

Let $X \subset \C^N$ be an affine 
complex analytic variety of dimension $n >1$, with an isolated 
complete intersection singularity at $0$. It is well known that the 
diffeomorphism type of its link $K = X \cap \sph_e$ does not depend  
on the choices of the embedding of $X$ in $\C^N$ nor on the sphere
$\sph_\e$, provided this is small enough. Moreover, according to \cite {Var}
one has a canonical contact structure $\cal X$ on
$K$, which is again independent of the embedding of $X$ in $\C^N$ and the 
choice of the sphere, up to contact-isomorphism. We refer to  $\cal X$ as the
{\it canonical contact structure} on $K$. Our interest in looking at this 
contact structure comes from \cite {CNP}, where the authors use it to make
 interesting applications to the theory of surface singularities.

This contact structure corresponds to the complex
sub-bundle $\D_X$ of the tangent bundle TK whose fiber at each point is the
complex orthogonal complement in $T(X \setminus \{0\})$ of the
 unit outwards normal vector field $\tau$ of $K$ in $X$.

For instance, if $n =2$ then  one has 
a  nowhere-vanishing holomorphic 2-form $\Omega$ around $0$ in $X$, 
which determines an $Sp(1)$-structure on the complex bundle $T(X-\{0\})$ 
(see \cite {Se1}).
If, as before, we denote by $\tau$ the unit outwards normal vector field of $K$ 
in $X$, then the bundle  $\D_X$ is the trivial 1-dimensional 
complex bundle spanned by the vector field
$j \cdot \tau$, obtained by multiplying the vector $\tau(x)$ by the quaternion 
$j$ at
each point of $K$. The vector field $i \cdot \tau$ is, up to scaling, the Reeb 
vector
field of the canonical contact structure $\cal X$.

 In this work we give a necessary and sufficient condition for $\D_X$
to be trivial  when $n > 2$:

\vv
{\it {\bf Theorem 1.} The complex bundle $\D_X$ that defines the canonical 
contact structure on $K$ is $C^\infty$  trivial as a 
complex vector bundle iff the Milnor number $\mu(X,0)$ 
of the germ  $(X,0)$ satisfies:
$$ \qquad \qquad \mu(X,0) \equiv \; (-1)^{n-1} \qquad mod\,(n-1)!\;.$$
 }

So, for instance, for the quadric $X = \{z_0^2 + \cdots + z_n^2 = 0\}$ the bundle
$\D_X$ is trivial if and only if  $n=2$ or $n$ is an odd number. 

Theorem 1 is a consequence of Theorem 2 below, which is 
proved via classical obstruction theory. We recall
that the GSV-index of a vector field $v$ on $X$, singular only at $0$, equals the 
Poincar\'e-Hopf index of a continuous extension of $v$ to a 
Milnor fibre of $X$ at $0$ (see \cite {GSV, BLSS}).

\vv
{\it {\bf Theorem 2.} 
Let $v$ be a 
nowhere-zero, continuous vector field on a neighbourhood of $M$ in $W$. Then 
the complex orthogonal complement of $v$ in $T(X \setminus \{0\})$ is a
$C^\infty$ trivial 
complex bundle if and only if the GSV-index of $v$ is a multiple of $(n-1)!$.}

\vv

Theorem 2 was announced in \cite {Se2}, with an outline of its proof. Here we give 
a  complete, self-contained proof of this result. This is also very much indebted 
to \cite {LS}, where similar arguments are used in relation with framed cobordism.

 We work always in the category of topological spaces and 
continuous maps, so our proofs of theorems 1 and 2 actually discuss topological 
triviality of the vector bundles in question. But everything becomes automatically 
$C^\infty$ because every continuous map between smooth manifolds can be 
approximated by a smooth map.

A holomorphic vector field $v$ on $X$, singular only at $0$ defines a 1-dimensional
holomorphic foliation on $X^* = X \setminus \{0\}$. In various situations one is 
naturally lead to considering the normal bundle of the foliation, and sections of 
it. Theorem 2 above implies the following:

\vv
{\it {\bf Corollary 1.}} Assume $X$ has complex dimension 3.
  Let $\cal F$ be a holomorphic,  locally free, 1-dimensional 
foliation on $X$, singular only at $0$, and let 
$\nu(\cal F)$ be its normal bundle in $X^*$. Let 
$v$ be a vector field on $X$  tangent to $\cal F$ and singular only at $0$.
 The following 
conditions are equivalent:

i) $\;\,$ The GSV-index of $v$ at $0$ is even.

ii) $\;$  The bundle $\nu(\cal F)$ admits a nowhere-zero $C^\infty$ section.

iii) $\,$The bundle $\nu(\cal F)$ is $C^\infty$ trivial.

\vv
When $X$ is smooth at $0$, the condition of $\cal F$ being locally free is
always satisfied, so it can be dropped. For $X$ smooth, the 
GSV-index is the usual local Poincar\'e-Hopf index of $v$.

\section{On compact parallelizable manifolds}

Let $W$ be a $2n$-dimensional, $n \ge 1$, compact, connected 
manifold with non-empty boundary $M$ and
 trivial tangent bundle $TW$. Let ${\cal F}_o$ be a trivialization of $TW$ 
and use this to define a complex structure on $TW$ and an isomorphism 
$TW \cong W \times \C^n$. Let $v$ be a 
nowhere-zero, continuous vector field on a neighbourhood of $M$ in $W$. It is 
well-known that $v$ can be extended to the interior of $W$ with isolated 
singularities, and the total number of singularities, 
counted with their 
corresponding local  Poincar\'e-Hopf indices, is independent 
of the choice of the extension of $v$ to the interior of $W$. This number is
the total Poincar\'e-Hopf index $\ind(v;W)$ of $v$ in $W$.

\begin{lemma}\label{theorem 1.1}
 If $\ind(v;W)$ is a multiple of $(n-1)!$, then $v$ can be 
completed to a continuous trivialization of the complex vector bundle $TW|_M$. 
That is, there exist  $(n-1)$ continuous sections 
$\a_2,...,\a_n$ of $TW|_M$, such that the set $\{v, \a_2,...,\a_n\}$ defines
a trivialization of 
$TW|_M$.
\end{lemma}

{\bf Proof.} Since $W$ is parallelizable and has non-empty boundary, there is 
 an immersion ${\cal I}: W \to {\R}^{2n}$, by
the immersion theorem of Hirsch-Poenar\'u \cite {Po}. Thus one has an induced 
(Gauss-type) continuous map 
$M \buildrel {\psi_v} \over {\longrightarrow} \sph^{2n-1}$, defined by
$$ \psi_v(x) \;=\, \frac{D {\cal I} (v(x))}{\vert D {\cal  I} (v(x)) \vert}\;,$$
where $D$ is the derivative.

By obstruction theory (see \cite {St}),
 $v$ can be extended to all of $W$ minus one point, say $x_o$, around which 
$\cal I$ can be assumed to be an embedding. Hence the topological degree 
of $\psi_v$ equals $\ind(v;W)$. In particular one has that some other vector field
$v'$ on a neighborhood of $M$ in $W$ is homotopic to $v$ 
(through never-vanishing vector fields) if and only if $\ind(v';W) = 
\ind(v;W)$. 

Now assume that $\ind(v;W)$ is a multiple of $(n-1)!$, {\it i.e.}, 
$\ind(v;W)= t (n-1)!$ for some integer $t$.
Notice that one has on $W$ vector fields with all possible Poincar\'e-Hopf
total indices and never-zero on $M$. Let $v_o$ be such a vector field with index
$t$. 

We recall (see \cite {Bo}) that the homotopy group $\pi_{2n-1}(U(n))$ of the 
unitary group $U(n)$ is 
isomorphic to $\Z$, and it has a canonical generator that maps to the generator 
$1$ of $\Z$. Let $\xi: \sph^{2n-1} \to U(n)$ represent this generator and define a map 
$$\phi_{v_o} : M \to U(n) \;,$$
by the composition $\phi_{v_o} \,=\, \xi \circ \psi_v$. Now, following 
\cite {Ke, KM},  twist the trivialization  ${\cal F}_o$ on the boundary $M$ 
using the map $\phi_{v_o}$; we get a new trivialization $\cal F$ of $TW|_M$. 
This means that at each point $x \in M$ we change the basis of $T_xW$given by
 ${\cal F}_o$ into its image by the linear map $\phi_{v_o}(x) \in U(n)$.
We claim that  $\cal F$ has $v$ as one of its $n$ sections, up to homotopy; 
this will complete the proof of the lemma.

To prove the above claim notice first that one has a map
$M \buildrel {\psi_{v_o}} \over {\longrightarrow} \sph^{2n-1}$ defined 
similarly to $\psi_v$ but now using $v_o$. The previous discussion implies 
that $\psi_{v_o}$ has degree $t$.

There is a fibration
$$U(n-1) \hookrightarrow U(n) \longrightarrow \sph^{2n-1}\,,$$
and an associated long exact homotopy sequence,
\begin{equationth}\label{exact seq}
 \cdots \to  \pi_{2n-1}(U(n)) \buildrel {p_*} \over {\longrightarrow} 
\pi_{2n-1}(\sph^{2n-1}) \to  \pi_{2n-2}(U(n-1)) \to \pi_{2n-2}(U(n)) \to 
\cdots
\end{equationth}

We know that $\pi_{2n-1}(\sph^{2n-1})  \cong \Z$.
Bott's calculations in \cite {Bo} tell us that:

i)  $\pi_{2n-1}(U(n)) \cong \Z$,

ii) $\pi_{2n-2}(U(n-1)) \cong \Z/(n-1)!$, 

iii) $\pi_{2n-2}(U(n)) \cong 0$ and
$p_*$ is multiplication by $(n-1)!$. 

Thus each section of $\cal F$ has index $[t \cdot (n-1)!]$ 
and the result follows. \qed

\section{On highly connected manifolds}

Now we assume that the manifold $W^{2n}$ of \S1  has the homotopy type of a 
bouquet of $n$-spheres, $n >1$. 
In this case one has the converse of \ref{theorem 1.1}:

\begin{lemma}\label{converse}
Let $v$ be a continuous section of $TW|_M$ which can be completed to a 
trivialization of the complex bundle  $TW|_M$; {\it i.e.}, there exist 
continuous sections $\a_2,...,\a_n$ of $TW|_M$ such that the set
$\F = \{v, \a_2,...,\a_n\}$ defines a trivialization of $TW|_M$ as a complex vector
 bundle. Then $\ind(v;W)$ is a multiple of $(n-1)!$. 
\end{lemma}
 
{\bf Proof.} We equip $W$ with
 a triangulation compatible with the boundary $M$, 
and we refer to $\F$ as {\it a complex framing} on $M$, meaning by this a 
trivialization of the complex bundle  $TW|_M$. We try to extend $\F$ to the 
interior of $W$ using the usual ``stepwise'' process: first to the $0$-skeleton,
then the $1$-skeleton and so on, as far as we can.

According to \cite {St}, the obstructions to extending $\F$ as a complex framing
over the interior of $W$ are elements in the relative cohomology 
$H^*(W,M; \Z)$. In fact these obstructions are all cocycles that represent,
 by definition,
the Chern classes of $W$ relative to the framing $\F$ on $M = \partial W$. 
Thus they live in the even-dimensional relative cohomology of $(W,M)$. 

By Lefschetz duality one has $H^i(W,M) \cong H_{2n-i}(W)$,
 hence all these groups vanish, except for $i = n, 2n$, since $W$ is assumed to
have the homotopy of a bouquet of $n$-spheres.

Let us assume first that $n$ is odd. Since Chern classes live in even dimensions, 
in this case the only possible obstruction to 
extending $\F$ to the interior of $W$ is the top relative Chern class 
$c_n(W,\F) \in H^{2n}(W,M;\Z)$. By definition, this class is the obstruction to 
extending to the interior of $W$ one of the sections that define $\F$, that we 
can take to be $v$. Hence $\F$ can be extended to all of $W$ minus one point, 
say $x_o$, and $\ind(v;W)$ can be regarded as being both, the local 
Poincar\'e-Hopf index at $x_o$ of the extension of $v$ to $W \setminus \{x_o\}$,
and also the Lefschetz dual  $c_n(W,\F)[W,M] \in H_0(W)$ 
of the Chern class $c_n(W,\F)$, where $[W,M]$ is the fundamental cycle of the pair.

Since
$\F$ is already extended to a trivialization of $T(W \setminus \{x_o\})$, 
one has that  $c_n(W,\F)$ can be identified, by excision, with the Chern class
of a small disc $D_\e$ in $W$ centered at $x_o$, relative to the framing
$\F$ on $\partial D_\e$. Then the exact sequence \ref{theorem 1.1} implies that
$c_n(W,\F)[(W,M)]$ is a multiple of $(n-1)!$, proving the lemma when $n$ is odd.

Consider now the case $n$ is even, say $n = 2m$, so $W$ has real dimension $4m$. 
We need: 

\begin{lemma}\label{lemma torus}
The framing $\F$ on $M$ extends to a trivialization $\widehat \F$
of the complex bundle 
$T(W\setminus \T)$, where $\T$ is a torus $\sph^n \times \sph^{n-1}$ embedded 
in the interior  $\buildrel {\circ} \over {W}$ of $W$.
\end{lemma}

{\bf Proof.}
Now there is a first possibly non-zero 
obstruction $\ro(\F) \in H^n(W,M;\Z)$ for extending $\F$ to the interior of $W$.
Let ${\cal S_{\F}} \in {H_{n}}(W)$ be the lefschetz dual of $\ro(\F)$. 

By hypothesis $W$ is simply connected. Hence Lemma 6 in \cite {Mi1} implies that
 $\cal S_{\F}$ can be represented by an $n$-sphere $S_{\F} $ 
embedded in  $\buildrel {\circ} \over {W}$. We claim:

i)  that $S_{\F}$ is actually embedded in  $\buildrel {\circ} \over {W}$
with trivial normal bundle, so the boundary of a tubular neighbourhood of 
$S_{\F}$ is homeomorphic to a torus $\sph^n \times \sph^{n-1}$; and

ii) the framing
$\F$ extends to a complex framing $\widehat \F$
on all of $W \setminus S_{\F}$.

\v

These two claims obviously prove the lemma 
\ref{lemma torus}. We prove (ii) first, and then use (ii) to prove that
$S_{\F}$  is embedded with trivial normal bundle.
We need:

\begin{lemma} The cohomology groups $H^i(W \setminus S, M)$ vanish for 
$ i = 0, 1, \cdots, n-2, n+2,\cdots,2n-2$; {\it i.e.}, for 
$i \ne n \pm 1, 2n-1,2n$.
\end{lemma}

Maybe this lemma can be improved, but this is all we need.

\v
{\bf Proof.} We will prove that $H^i(W \setminus S, M)$ and $H^i(W, M)$ are
isomorphic in the above dimensions. The lemma then follows from the fact
that $W$ is (homotopicaly) a bouquet of $n$-spheres.

Consider the exact cohomology sequences of the pairs $(W,M)$ and  
$(W \setminus S_{\F},M)$. The ``Five Lemma'' implies that
$H^i(W,M)$ is isomorphic to $H^i(W \setminus S_{\F} ,M)$ whenever the groups 
$H^j(W)$ and $H^j(W \setminus S_{\F})$ are isomorphic  for $j = i-1$ and $j= i$.
 Then the exact sequence of the pair $(W, W\setminus S_{\F})$, 
together with Alexander duality and the fact that $ S_{\F}$ is an $n$-sphere, 
imply the lemma. \qed

\vv

Let us return to the proof of (ii). Let $N$ be an open tubular neighbourhood of
 $ S_{\F}$ in $W$ and consider a triangulation $T$ of the compact manifold
$W \setminus N$, compatible with its boundary $M \cup N$. Denote by $T^{(i)}$ 
the corresponding $i^{th}$-skeleton. 

Now we start the process of trying to extend $\F$ from $M$  to all of 
$W \setminus N$, step by step.
 By the previous lemma, there are no obstructions up to dimension
$(n-2)$, and $(n-1)$ is odd, so we can extend $\F$ to $M \cup T^{({n-1})}$ 
because the obstructions  appear in even dimensions. 

Notice one has a commutative diagram:
\begin{equationth}\label{no-obstruction}
\xymatrix{
H^n(W, W \setminus S_{\F} ) \ar[d]^L \ar[r]^-{j^*} & H^n(W, M) 
\ar[d]^P \ar[r]^-{k^*} &
H^n(W \setminus S_{\F}, M)\\
H_n(S_{\F}) \ar[r]^{i_*} & H_n(W) & \\
}
\end{equationth}
where the horizontal arrows are induced by the inclusions,
 $L, P$ denote Lefschetz and Poincar\'e duality and the first row is exact in 
the middle.  This implies that the obstruction class 
$\ro(\F)$ is in the image of $j^*$, so it is mapped to $0$ by 
$k^*$, thus implying that there is no obstruction to extending 
$\F$  to a complex framing $\widehat \F$
on the $n$-skeleton of the triangulation on $W \setminus N$.

Now,  $n+1$ is odd, so we can extend the framing over $M \cup T^{n+1}$, and then 
the previous lemma grants we can continue the process up to dimension $2n-2$, 
and therefore to dimension $2n-1$ since this is odd.

We finally meet an obstruction in the top dimension $2n$.  Regarded in 
homology, this means the obstruction has support at isolated points, say
$q_1,\cdots,q_r$ which we can assume to be in the interior of $W \setminus N$. 
We can now retract $N$ to $S_{\F}$ smoothly, thus extending $\F$ to all
of $W \setminus (S_{\F} \cup \{q_1,\cdots,q_r\})$. Finally,
 we can move these points in $W$ by an isotopy
 so that they become contained in the 
sphere  $S_{\F}$, and we arrive to statement (ii).

%%%%%%%%%%%%%%%%%%%%%%%%%%%%%%%%%%%%%%

 To prove the first claim we notice that, since $W$ is parallelizable,
 Lemma 7 in \cite {Mi2} implies that 
$S_{\F}$ is embedded in $\buildrel {\circ} \over {W}$ with trivial normal 
bundle if and only if
its self-intersection number $S_{\F} \cdot S_{\F}$ is $0$, that is if and only 
if $$\ro(\F) \cup \ro(\F) \,=\,0\;.$$
Let us prove that this happens. Recall we have the trivialization 
 ${\cal F}_o$ of $TW$ that we used to define a complex structure on $TW$. 
Of course the Chern classes $c_i(W,M; {\cal F}_o)$
of $W$ relative to the framing ${\cal F}_o$ on 
$TW|_M$ are all zero. 

In general one has that given two complex framings 
$\b, \b'$ on the boundary $M$, each defines relative Chern classes 
$c_i(W;\b), c_i(W;\b')$ in $H^{2i}(W,M)$ for all $i \ge 0$ even, and in each 
dimension their difference is a contribution of the boundary
$$ H^{2i-1}(M) \buildrel {\delta^*} \over \longrightarrow H^{2i}(W,M)\;,$$
given by the {\it difference cocycle} (see \cite {St}, also \cite {Ke2,BLSS}).

In our case, since the classes $c_i(W,M; {\cal F}_o)$ are all zero, this implies
that $c_m(W,M; \cal F)$ is a coboundary, {\it i.e.}, $\ro(\F) = \delta^*(a)$ 
for some $a \in H^{2m-1}(M)$. Thus one has:
$$\ro(\F) \cup \ro(\F)\,=\,\delta^*(a) \cup \ro(\F) \,=\,
j^* \delta^*(a) \cup \ro(\F)\,=\,
0\;,$$
where $j^*: H^{*}(W,M) \to H^*(W)$ is induced by the inclusion, 
proving lemma \ref{lemma torus}. \qed

\v

Let us now complete the proof of lemma \ref{converse} when $n$ is even,  
$n =2m$. 

 From the previous lemma we know that $\widehat \F$ is a complex framing 
that extends $\F$ to all of $W$ minus the interior 
${\rm Int}\, {\widehat \T}$ 
of the solid torus 
$\widehat \T  \cong \sph^n \times \B^{n} $, which is a 
tubular neighbourhood of the  $n$-sphere $S_{\F}$ in $W$. Set 
$\T = \partial \widehat \T$.

We know, 
essentially by definition, that
$$c_n(W;\F)[W,M] = \ind (v;W)\,,$$
where $[W,M]$ is the fundamental cycle of the pair (W,M), 
and 
$$\,c_n(W;\F)[W,M] = c_n(\widehat \T, \widehat \F)[\widehat \T, \T]$$
because $\widehat \F$ extends $\F$. Let us prove that the latter integer
is a multiple of $(n-1)!$. 

 Let ${\cal F}_o$ be the trivialization of $TW$ used to determine an 
almost-complex structure on $W$. Since $\pi_{n}(U( n)) = 0$ (see \cite {Bo}), 
$n > 1$, 
we can assume that  ${\cal F}_o$ and  $\widehat \F$ coincide, up to homotopy, 
over a parallel $(\sph^n \times *)$ of $\T$, 
where $*$ is a point in $\partial \B^{n}$.

Recall that, following \cite {Ke}, we can think of the framing 
 $\widehat \F|_{\T}$ as being obtained from  ${\cal F}_o|_\T$, twisting 
it by a map ${\mathfrak  f}: \T \to U(n )$. This map is, at each point $z$ of $\T$, 
the linear transformation that carries the  base of $T_zW$ given by  
${\cal F}_o(z)$ into the one given by $\widehat \F(z)$. In the sequel, we 
identify the framing $\widehat \F|_{\T}$ with the map  ${\mathfrak  f}$.

Notice that $\T \cong \sph^n \times \sph^{n-1}$ has a CW-decomposition as 
$e^0 \cup e^{n-1} \cup e^n \cup e^{2n}$,
 with $\sph^n = e^0 \cup  e^n$ and
 $\sph^{n-1}=  e^0 \cup e^{ n-1}$. 

Now define a 
 complex framing $\b$ on $\T$ by 
$$\b(x,y) = \widehat \F (x_o,y)\,,$$
where $x_o$ is the projection of $e^o$ to  $\sph^n$ and we are thinking of 
 $\widehat \F|_{\T}$ as being the map  ${\mathfrak  f}$. 
Then $\b$ and $\widehat \F|_{\T}$ agree over $e^0 \cup e^{n-1} \cup e^n$, 
by construction, and therefore they differ by a map 
$\partial e^{2n} \to U(n )$, {\it i.e.}, by an element in $\pi_{2n-1}(U(n))$.

We claim that one has 
$c_n(\widehat \T;\b)[\widehat \T,\T] = 0$. For this, it is convenient to think 
of this integer as being the degree of $\b$, {\it i.e.}, if we write $\b$ as
the $n$-frame $\{\xi_1,...,\xi_n\}$, where the $\xi_i$ are vector fields, then
 $c_n(\widehat \T;\b)[\widehat \T,\T] = \ind(b_i;\widehat \T)$ for all $i =1,...n$.
Take one of these vector fields, say $b_1$. Notice that the framing 
${\cal F}_o$ provides an isomorphism $TW \cong W \times \C^n$. Restricting this to
$\widehat \T$ and projecting into $\C^n$ one gets a well defined map 
$\hat \b = p \circ b_1$,
$$\T \buildrel {\frac{b_1}{|b_1|}} \over \longrightarrow T \widehat \T 
\longrightarrow \sph^{2n-1}\;,$$
whose degree equals $\ind(b_i;\widehat \T)$. To prove that this map has degree 
$0$ it 
is enough to show that $\hat \b$ extends to a map $ \widehat \T \to \sph^{2n-1}$ 
(see for instance \cite {GP}). This follows because $\b$ is constant on
$\sph^n \times \{*\}$ and  $\hat \b|_{e^0 \cup e^{n-1}}$ is nulhomotopic because
$\sph^{2n-1}$ is $(2n-2)$-connected.

Therefore we have that $c_n(\widehat \T;\b)[\widehat \T,\T] = 0$ and, 
by construction, $\b$ differs from $\widehat \F$ on $\T$ by an element
in $\pi_{2n-1}(U(n))$. The result now follows from the exact sequence 
\ref{exact seq}, where $p^*$ is multiplication by $(n-1)!$. \qed

\section{Proof of the theorems}

We now take the manifold $W$ to be a Milnor fiber of $X$. More precisely, we assume
the $n$-dimensional complete intersection 
germ $(X,0)$ is defined by a reduced function 
$$f = (f_1,\cdots,f_k): (U,0) \to (\C,0)\;,$$
with a critical point at $0$,
where $U$ is an open neighbourhood of $0$ in $\C^{n+k}$, and we let
$W = f^{-1}(t) \cap \B_\e$, where $\B_\e$ is a ball around $0 \in \C^{n+k}$ of 
sufficiently small radius $\e >0$, and $t$ is a regular value with $|t|$ 
sufficiently small with respect to $\e$. We are assuming also that $0$ is an 
isolated singularity in $X$, which means it is an isolated critical point of
$f$ in $X$.

We know from \cite {Mi2, Ha} (see also \cite {Lo}) that $X$ is parallelizable
and it has the homotopy type of a bouquet of $\mu = \mu(X,0)$ spheres of dimension
$n$, for some integer $\mu > 0$ which is known as {\it the Milnor number} of
$(X,0)$. In fact one has more: the tangent bundle $TW$ is naturally equipped with a
 complex structure, since $W$ is defined by a holomorphic function, and 
one has:
$$T\C^{n+k}|_W \,\cong \, TW \oplus \nu W\;,$$
as $C^{\infty}$ vector bundles, where $\nu W$ is the normal bundle of $W$. 
The bundle $\nu W$ is 
canonically trivialized by the gradient vector fields 
$(\nabla f_1,\cdots,\nabla f_k)$, where
$$\nabla f_i = \Big(\frac{\overline {\partial f_i}}{\partial z_1},\cdots, 
\frac{\overline {\partial f_i}}{\partial z_{n+k}}\Big)\,.$$
This implies that $TW$ is {\it stably trivial} as a complex vector bundle 
(see \cite {KM}), {\it i.e.}, its tangent bundle plus a trivial bundle is trivial. 
Thus, by \cite {KM}, it is actually trivial as a complex bundle,
because $W$ is connected with non-empty boundary and therefore $H^{2n}(W) \cong 0$.

Thus we are in the situation envisaged in sections 1 and 2 above, and Theorem 2
follows from lemmas \ref{theorem 1.1} and \ref{converse}, together with the fact
that the Thom isotopy theorem (see for instance \cite {AR}) allows
us to identify the link $K = X \cap \partial \B_\e$ with the boundary of 
 $W$ (see \cite {Mi2, Lo}).  \qed

\vv
Now take $v$ to be a radial, outwards-pointing vector field on $X$. Then its 
GSV-index equals its Poincar\'e-Hopf index in $W$, where we are identifying
$K$ with the boundary of $W$. In this case $v$ is everywhere transversal to
$\partial W$, pointing outwards. Hence,   by the theorem
of Poincar\'e-Hopf for manifolds with boundary, its index in $W$ is $\chi(W)$,
the Euler-Poincar\'e characteristic. By \cite {Mi2, Ha} 
(see also \cite {Lo}),  this 
equals $ 1 + (-1)^{n} \mu(X,0)$, since $W$ has the homotopy type of a bouquet 
of $ \mu(X,0)$ spheres of dimension $n$. 
Hence Theorem 2 implies
$$1 + (-1)^{n}\mu(X,0) \equiv 0 \qquad \hbox{mod} \; (n-1)!\;,$$
and we arrive to Theorem 1.  \qed

\section{Proof of Corollary 1 and an example}

Now we prove Corollary 1. 
Notice first that $X$ is locally a cone over its link $K$ with vertex at $0$, 
by \cite {Mi2}. Hence in statements (ii) and (iii) of Theorem 3 we can restrict the 
discussion to the link.

The equivalence between statements (i) and (iii) is immediate from Theorem 2, and 
it is obvious that (iii) implies (ii), so we only must prove that (ii) implies 
(iii). Let $\xi$ be a  never-zero continuous
 section of the normal bundle $\nu(\cal F)$. This spans a 1-dimensional 
continuous complex line sub-bundle $\cal L$ of $\nu(\cal F)$. 
The bundle  $\cal L$ is trivial iff  $\nu(\cal F)$ is trivial. But $K$ is 
2-connected, by \cite {Mi2}. Hence every complex line bundle over $K$ is trivial.

\vv

To finish this note, we give an example. Let $X$ be a hypersurface in $\C^{2n}$
defined by some function $f: (\C^{2n},0) \to (\C,0)$. Then the Hamiltonian vector 
field 
$$v = \Big( \frac{\partial f}{\partial z_2}\,,\, -\frac{\partial f}{\partial z_1}
\,, \, \frac{\partial f}{\partial z_4 }\, ,\, -\frac{\partial f}{\partial z_3 }\,,\,
 \cdots\,,  \frac{\partial f}{\partial z_{2n}}\,,\,
 -\frac{\partial f}{\partial z_{2n-1}} \Big)
$$
is obviously tangent to $X^*$, since $df(v) \equiv 0$ everywhere. This implies 
also that $v$ is tangent to all fibers of $f$. Hence its GSV-index is $0$.
Thus, by Theorem 2, the normal bundle of the holomorphic foliation that $v$ spans
on  $X^*$ is topologically trivial.

\bibliographystyle{plain}

$\,$

Jos\'e Seade

Instituto de Matem\'aticas, UNAM, 

Unidad Cuernavaca,

A. P. 273-3, Cuernavaca, Morelos, 

M\'exico.

jseade@matcuer.unam.mx

\end{document}